\documentclass[12pt,titlepage]{article}

\usepackage[tbtags]{amsmath}
\usepackage{amsopn,amsthm,amsfonts}

\DeclareMathOperator{\Res}{Res}

\newtheorem{theorem}{Theorem}

\theoremstyle{definition}

\theoremstyle{remark}

\newtheorem*{conjecture}{Conjecture}
\newtheorem*{acknow}{Acknowledgments}

\begin{document}

 \title{A Rodrigues formula for the Jack polynomials and the Macdonald-Stanley
conjecture}
\author{Luc~\textsc{Lapointe} and Luc~\textsc{Vinet}\\
\small
\begin{tabular}{c}
Centre de recherches math{\'e}matiques\\
Universit{\'e} de Montr{\'e}al\\
C. P.~6128,  succ.~Centre-ville\\
Montr{\'e}al QC H3C 3J7
\end{tabular}}
\date{\textbf{CRM-2294}\\[\bigskipamount]
  July 1995}
\maketitle

\begin{abstract}
  A formula of Rodrigues-type for the Jack polynomials is presented.
  It is seen to imply a weak form of a conjecture of Macdonald and
  Stanley.
\end{abstract}

\section{Introduction}

As usual, a partition $\lambda$ of $N$ will be a sequence $\lambda =
(\lambda_1, \lambda_2, \dots )$ of non-negative integers in decreasing order
$\lambda_1 \ge \lambda_2 \ge \dots$ such that $|\lambda| = \lambda_1 +
\lambda_2 + \dots = N$.  The integers $\lambda_i$ are called the parts of
$\lambda$.  Let $\lambda$ and $\mu $ be two partitions of $N$.  In the
dominance ordering, $\lambda \ge \mu $ if $\lambda_1 + \lambda_2 + \dots +
\lambda_i \ge \mu _1 + \mu _2 + \dots + \mu _i$ for all $i$.

Let $\Lambda_n$ denote the ring of symmetric functions in the variables $x_1,
x_2, \dots , x_n$.  Two natural bases of $\Lambda_n$ are: the power sums
$p_\lambda = p_{\lambda_1} p_{\lambda_2} \dotsm $ where $p_i = \sum_{k} x_k^i$
and the symmetric monomials $m_\lambda = \sum x_1^{\lambda_1} x_2^{\lambda_2}
\!\dotsm$ summed over distinct permutations.

To the partition $\lambda$ with $m_i(\lambda)$ parts equal to $i$, we associate
the number
\begin{equation} \label{1}
z_\lambda
        = 1^{m_1} m_1 !  \, 2^{m_2} m_2! \dotsm
\end{equation}
Let $\alpha$ be a parameter and $\mathbb Q(\alpha)$ the field of all rational
functions of $\alpha$ with rational coefficients.  A scalar product $\langle \
, \ \rangle$ on $\Lambda_n \otimes \mathbb Q(\alpha)$ is defined by
\begin{equation}
\langle p_\lambda, p_\mu \rangle
        =\delta_{\lambda \mu } z_\lambda \alpha^{\ell(\lambda)},
\end{equation}
where $\ell(\lambda)$ is the number of parts of $\lambda$.  The Jack
polynomials $J_\lambda (x_1, \dots , x_n; \alpha) \in \Lambda_n \otimes \mathbb
Q(\alpha)$ are uniquely specified by \cite{1,2}
\begin{align}
\mathrm{(i)} \ &  \langle J_\lambda, J_\mu \rangle = 0, \qquad \text{if }
\lambda \ne \mu , \label{3}\\
\mathrm{(ii)} \ &  J_\lambda = \sum_{\mu \le \lambda} v_{\lambda\mu }(\alpha)
m_\mu ,  \label{4}\\
\mathrm{(iii)} \ &  \text{if } |\lambda| = N, \quad v_{\lambda, 1^N} = N!
\label{5}
\end{align}
The Jack polynomials $J_\lambda(x_1, \dots , x_n;\alpha)$ where $\ell(\lambda)
\le n$ are also known to obey the differential equations \cite{3,1}
\begin{equation}\label{6}
H(\alpha) J_\lambda (x; \alpha)
        = \varepsilon_\lambda (\alpha) J_\lambda(x;\alpha),
\end{equation}
where
\begin{equation} \label{7}
H(\alpha)
        = \alpha \sum_{j=1}^n  \biggl( x_j \frac{\partial}{\partial x_j}
\biggr)^2 + \sum_{j < k} \biggl( \frac{x_j + x_k}{x_j - x_k} \biggr) \biggl(
x_j  \frac{\partial}{\partial x_j} - x_k \frac{\partial}{\partial x_k} \biggr),
\end{equation}
and
\begin{equation} \label{8}
\varepsilon_\lambda (\alpha)
        = \sum_{j=1}^n \bigl[ \alpha \lambda_j^2 + (n+1-2j)\lambda_j \bigr].
\end{equation}
The polynomials $J_\lambda$ appear in the wave functions of the
Calogero-Sutherland model \cite{4}.  This is an exactly solvable
quantum-mechanical system which describes $n$ particles on a circle interacting
pairwise through long-range potentials.  It is nowadays quite useful in studies
of phenomena associated with fractional statistics and there is considerable
interest in identifying the algebraic structure underlying this model.  In this
connection, it is natural to look for a formula giving the wave functions of
the excited states through the action of creation operators on the ground state
wave function.  We obtained \cite{3} as a result a formula of Rodrigues-type
for the Jack polynomials which implies a weak form of a longstanding conjecture
due to Macdonald and Stanley, namely that the coefficients $v_{\lambda\mu
}(\alpha)$ in \eqref{4} are polynomials in $\alpha$ with integer coefficients.
We present and discuss this formula here.  It will be proved elsewhere
\cite{3}.

\section{Creation operators}
The creation operators are constructed from the Dunkl operators \cite{5}
\begin{equation} \label{9}
D_i
        = \alpha x_i \frac{\partial}{\partial x_i} +
\sum_{\begin{subarray}{l}j=1\\ j\ne i \end{subarray}}^n \frac{x_i}{x_i - x_j}
(1 - K_{ij}), \qquad i = 1,2, \dots , n,
\end{equation}
where $K_{ij} = K_{ji}$ is the operator that permutes the variables $x_i$ and
$x_j$:
\begin{equation} \label{10}
K_{ij}x_i
        = x_j K_{ij}, \quad K_{ij} D_i = D_j K_{ij}, \quad K_{ij}^2 = 1.
\end{equation}
Let $J = \{  j_1,  j_2, \dots , j_\ell \}$ be sets of cardinality $|J| = \ell$
made of integers $j_\kappa \in \{ 1, \dots , n \}$, $1 \le \kappa \le \ell$
such that $j_1 < j_2 < \dots < j_\ell$ and introduce the operators
\begin{equation} \label{11}
D_J
        = (D_{j_1} + 1) (D_{j_2} + 2) \dotsm (D_{j_\ell} + \ell),
\end{equation}
labelled by such sets.  The creation operators $B_i^+$ are defined by
\begin{equation} \label{12}
B_i^+
        = \sum_{\begin{subarray}{c} J \subset \{ 1,\dots , n\}\\ |J| = i
\end{subarray}}  x_J D_J,
\end{equation}
with
\begin{equation} \label{13}
x_J
        = \prod_{i \in J} x_i.
\end{equation}
The sum in \eqref{12} is over all subsets $J$ of $\{ 1, \dots , n \}$ that are
of cardinality $i$.

\section{The Rodrigues formula}
We can now state our main result \cite{3}.

\begin{theorem} The Jack polynomials $J_\lambda(x;\alpha)$ associated to
partitions $\lambda = (\lambda_1 , \dots , \lambda_n)$ are given by
\begin{equation} \label{14}
J_\lambda(x;\alpha)
        = (B_n^+)^{\lambda_n} (B_{n-1}^+)^{\lambda_{n-1} -\lambda_n} \dots
(B_1^+)^{\lambda_1 - \lambda_2} \cdot 1.
\end{equation}
\end{theorem}

That the right-hand side of \eqref{14} is a symmetric polynomial of degree $N =
\lambda_1 + \dots + \lambda_n$ in the variables $x_1, \dots , x_n$ is easily
seen from the properties of the Dunkl operators. It is checked that the
operators \eqref{9} satisfy the commutation relations
\begin{equation} \label{15}
[D_i, D_j]
        = (D_j - D_i)K_{ij}.
\end{equation}
Let us now use the notation $\Res^{\{i,j,k, \dots\}} X$ or $\Res^J X$ to
indicate that the operator $X$  is taken to act on functions that are symmetric
in the variables $x_i, x_j, x_k, \dots$ or $x_{j_1}, x_{j_2}, \dots , j_\kappa
\in J$, respectively.  With $m$ some integer, it is straightforward to verify
with the help of \eqref{15} that
\begin{equation} \label{16}
\Res^{\{i,j\}}(D_i + m)(D_j + m+1)
        =\Res^{\{i,j\}} (D_j+m)(D_i+m+1).
\end{equation}
It follows that $\Res^JD_J$ is invariant under the permutations of the
variables $x_{j_\kappa}$, $j_\kappa \in J$ and that this operator therefore
leaves invariant the space of symmetric functions in these variables.

Recalling how the creation operators $B_i^+$ are constructed in terms of the
operators $D_J$, it is then clear that $(B_n^+)^{\lambda_n} \!\dots
(B_1^+)^{\lambda_1 - \lambda_2} \cdot 1$ is a symmetric function of the
variables $x_1, \dots , x_n$.  That it is a homogeneous polynomial of degree
$N$ is readily seen from observing that $B_i^+$ has scaling dimension $i$.

Let $\varphi_{(\lambda_1,\dots ,\lambda_i,0,\dots )}
        = (B_i^+)^{\lambda_i} \dots  (B_1^+)^{\lambda_1 - \lambda_2} \cdot 1$.
The proof of the theorem involves showing that
\begin{equation} \label{17}
[H(\alpha),B_i^+] \varphi_{(\lambda_1,\dots ,\lambda_i,0,\dots )}
        = B_i^+ \biggl\{ 2\alpha \sum_{j=1}^n \lambda_j + i \alpha + i(n-i)
\biggr\}
\varphi_{(\lambda_1, \dots, \lambda_i,0,\dots )}.
\end{equation}
If \eqref{17} is true, it is clear that successive applications of the creation
operators on 1 will build eigenfunctions of $H$.  One then iteratively checks
that the spectrum coincides with \eqref{8} to confirm the identification
\eqref{14}.  A lengthy combinatorial argument which will be published elsewhere
\cite{3} actually leads to \eqref{17}.

\section{The conjecture of Macdonald and Stanley}
Many conjectures involving the Jack polynomials have been formulated.  A famous
one is due to Macdonald and is reproduced in Stanley's reference
article.\cite{1}  It is stated as follows.

\begin{conjecture} Let
\begin{equation} \label{18}
\widetilde v_{\lambda \mu }(\alpha)
        = \frac{v_{\lambda\mu }(\alpha)}{\prod_{i \ge 1} m_i(\mu )!},
\end{equation}
where $v_{\lambda\mu }(\alpha)$ are as in \eqref{4}.  Then $\widetilde
v_{\lambda \mu }(\alpha)$ are polynomials in $\alpha$ with nonnegative integer
coefficients.
\end{conjecture}

To our knowledge, it was still not even known whether $v_{\lambda \mu
}(\alpha)$ are polynomials.  This can now be proved and the following weak form
of the above conjecture is seen to follow directly from fomula \eqref{14}.

\begin{theorem}
The coefficients $v_{\lambda\mu }(\alpha)$ in \eqref{4} are polynomials in
$\alpha$ with integer coefficients.
\end{theorem}

The proof proceeds in a recursive fashion.  Assume that the assertion is true
for partitions $\lambda = (\lambda_1, \dots , \lambda_i, 0, \dots)$.  Formula
\eqref{14} allows to increase the values of the parts and their number by
acting on $J_\lambda(x;\alpha)$ with $B_k^+$, $k \ge i$.  Theorem~2 will be
proved in these higher cases and thus in general since one can start with
$J_0(x;\alpha) = 1$, if the matrix elements of the creation operators in the
symmetric monomial basis are shown to be polynomials in $\alpha$ with integer
coefficients.  To convince oneself of that, let $m_\mu = \sum_{\text{distinct
perm}} \widehat m_\mu $ with $\widehat m_\mu = x_1^{\mu _1} \dotsm   x_n^{\mu
_n}$ be one such symmetric monomial.  We know that $B_k^+ m_\mu =
\sum_{\text{distinct perm}} B_k^+ \widehat m_\mu $ is a symmetric function.
Moreover, the operators $D_m$, $m = 1, \dots , n$, and hence the operators
$B_k^+$ are easily found to give when acting on monomials $\widehat m_\mu $,
sums of similar monomials multiplied

 by polynomials in $\alpha$ with integer coefficients.  Owing to the symmetry
of $B_k^+ m_\mu $, $\sum_{\text{distinct perm}} B_k^+ \widehat m_\mu $ must
produce the desired result for the matrix elements of the creation operators in
the symmetric monomial basis.

\section{Conclusion}
We believe that the Rodrigues formula that we have obtained will provide a
useful tool to further advance the proofs of the outstanding conjectures on the
Jack polynomials.  We also trust that it will help obtain the dynamical algebra
of the Calogero-Sutherland model.  We are currently developing the extensions
of the results presented here to the Macdonald polynomials \cite{2} as well as
to the multivariate special functions associated to lattices other than the one
associated to $A_{n-1}$.  We hope to report on these issues in the near future.

\begin{acknow}
We would like to express our thanks to Fran{\c{}}{c}ois Bergeron, Adriano
Garsia and Doron Zeilberger for various comments and suggestions.

\noindent This work has been supported in part through funds provided by NSERC
(Canada) and FCAR (Qu{\'e}bec).  L.~Lapointe holds a NSERC postgraduate
scholarship.
\end{acknow}


\begin{thebibliography}{33}
\bibitem{1}
R.~P. Stanley, \emph{Some combinatorial properties of Jack symmetric
functions}, Adv. Math. {\bf 77} (1988), 76--115.
\bibitem{2}
I.~G. Macdonald, \emph{Symmetric functions and Hall polynomials}, 2nd edition,
Clarendon Press, Oxford, 1995.
\bibitem{3} L.~Lapointe and L.~Vinet, \emph{Exact operator solution of the
Calogero-Sutherland model}, CRM-preprint {\bf 2272} (1995).
\bibitem{4} D.~Haldane, \emph{Physics of the ideal fermion gas: spinons and
quantum symmetries of the integrable Haldane-Shastry spin chain}, Correlation
Effects in Low-Dimensional Electron Systems (A. Okiji and N.~Kamakami, eds.),
Springer-Verlag, New York, 1995, pp. 3--20, and references therein.
\bibitem{5} C.~F. Dunkl, \emph{Differential-difference operators associated to
reflection groups}, Trans. Amer. Math. Soc. {\bf 311} (1989), 167--183.
\end{thebibliography}
\end{document}